\documentclass[12pt]{article}
\usepackage[english]{babel}
\usepackage{amssymb,amsmath,graphics}
\parskip\medskipamount
\newtheorem{theorem}{Theorem}[section]

\newtheorem{lemma}[theorem]{Lemma}
\newtheorem{cor}[theorem]{Corollary}
\def\<{\langle}
\def\>{\rangle}
\newcommand{\proof}{\emph{Proof.~}}
\newcommand{\cC}{\mathcal{C}}

\newcommand{\PG}{\mathsf{PG}}

\newcommand{\sfH}{\mathsf{H}}
\newcommand{\cS}{\mathcal{S}}

\newcommand{\cP}{\mathcal{P}}
\newcommand{\cG}{\mathcal{G}}
\newcommand{\cN}{\mathcal{N}}

\def\qed{{\hfill\hphantom{.}\nobreak\hfill$\Box$}}

\newcommand{\cB}{\mathcal{B}}

\newcommand{\K}{\mathbb{K}}

\begin{document}

\author{Fabienne Haot \and Koen Struyve\thanks{The second author is supported by the Fund for Scientific Research --
Flanders (FWO - Vlaanderen)} \and Hendrik Van Maldeghem\thanks{The third author is partly supported by a
Research Grant of the Fund for Scientific Research -- Flanders (FWO - Vlaanderen)} }
\title{\bf Ree Geometries}

\date{}
\maketitle

\begin{abstract}
We introduce a rank 3 geometry for any Ree group over a not necessarily perfect field and show that its full
collineation group is the automorphism group of the corresponding Ree group. A similar result holds for two rank 2 geometries obtained as a truncation of this rank 3 geometry. As an
application, we show that a polarity in any Moufang generalized hexagon is unambiguously determined by its set
of absolute points, or equivalently, its set of absolute lines.
\end{abstract}

{\footnotesize AMS Classification: 20B25, 51E05, 51E12}

\section{Introduction}

The Ree groups in characteristic 3 (defined by Ree in \cite{Ree:61}) and their generalizations over non-perfect
fields (by Tits \cite{Tit:60}) provide examples of rank 1 groups, or Moufang sets, or split BN-pairs of rank 1
that behave roughly as simple algebraic groups of rank 1, but over a one-dimensional root system consisting of
six roots, since the unipotent radicals have nilpotency class 3. This is a rather rare phenomenon; indeed,
until recently, these were the only known rank 1 groups with this property (a second class was discovered and
constructed in \cite{Mue-Mal:06}).  Associated with each Ree group is a geometry (called a \emph{unital} in the
finite case), which is a linear space (in the finite case a $2-(q^3+1,q+1,1)$-design), see \cite{Lun:66}. This
geometry can be viewed as the geometry of involutions in a Ree group, since the blocks are in one-to-one
correspondence with a conjugacy class of involutions (in the finite case there is only one conjugacy class).
This way, Ree groups can be better understood in that several properties become more geometric and intuitive
through this geometry.

In the present paper we introduce another geometry for each Ree group, inspired by the general construction of
geometries associated to ``wide'' rank 1 groups as proposed by Tits in one of his lectures; see \cite{Mal:07}
for an account on this idea (``wide'' here means that the unipotent subgroups are not abelian). In fact, this
construction is the counterpart for Ree groups of the inversive planes for Suzuki groups (see also
\cite{Mal:07}). The structure of the geometries that we will introduce is probably slightly more involved than
that of the ``unitals'', but they have the major advantage that the automorphism groups of the corresponding Ree groups are their full
automorphism groups (and this is our Main Result below), a result that is not yet proved for the unitals. This
result contributes to Tits' programme of characterizing all ``wide'' rank one groups in this way. As an
application, we can show that every collineation of a Moufang hexagon of mixed type permuting the absolute
points of a polarity, centralizes that polarity (or, equivalently, also permutes the absolute lines). This, in
turn, means that the set of absolute points of any polarity of any  Moufang hexagon (necessarily of mixed type)
determines the polarity completely and unambiguously.

The ``new'' geometries also have a number of interesting combinatorial properties, but we will not concentrate
on these, though it would be worthwhile to perform an investigation in that direction.

Every Ree group is the centralizer of a certain outer involution of a Dickson group of type $\mathsf{G}_2$ over
a field of characteristic 3 admitting a Tits endomorphism. A geometric way to see this is to consider the
associated Moufang generalized hexagon, which is of mixed type. Then the outer involution is a polarity, and the associated Ree group
acts doubly transitively on the absolute points of that polarity. That is exactly the way we are going to
define and use the Ree groups. These Moufang hexagons are called \emph{Ree hexagons} in \cite{Mal:98} precisely
for that reason.

Hence, in order to investigate the Moufang sets (or rank 1 groups) associated with the Ree groups, we turn to
the Ree hexagons, which, as follows from our remarks above, are defined over a field of characteristic 3
admitting a Tits endomorphism $\theta$, and they allow a polarity $\rho$. The absolute points under this
polarity, together with the automorphisms of the mixed hexagon commuting with $\rho$, form the \emph{Ree-Tits
Moufang set}.  Since we will need an explicit description of the absolute points of $\rho$, we will use
coordinates. These will be introduced in Section~\ref{hexagon}. We define the Ree geometries in
Section~\ref{Reegeometry} and state our main results and main application in Section~\ref{mainresult} (but we
formulate our main results also below in rough terms).  The rest of the paper is then devoted to the proofs.

Since the Ree groups have unipotent subgroups of nilpotency class 3 (at least, if the base field is large
enough), the Ree geometries that we will define have rank 3. This means that we will have two types of
\emph{blocks} in our geometry. In this paper we prove that every automorphism of such a geometry is an
automorphism of the corresponding Ree group, by writing down explicitly the automorphisms of this geometry. But
we also do slightly better and prove that the same conclusion holds when restricting to one type of blocks. We
call these geometries \emph{truncated Ree geometries}.  Hence, loosely speaking, we may write our main result
as follows:

\begin{quote}\emph{
The full automorphism group of a (truncated) Ree geometry is induced by the full collineation group of the
corresponding Ree hexagon.}
\end{quote}

\section{The coordinatization of the Ree hexagon}\label{hexagon}
In this section, we present two coordinatizations of the Ree hexagons, which can at the same time serve as a
definition of these structures. We start with the coordinatization with respect to one flag
$\{(\infty),[\infty]\} $ (a \emph{flag} is a pair of elements consisting of a point and a line that are
incident).  This coordinatization was first carried out by De Smet and Van Maldeghem for (finite) generalized
hexagons in \cite{Sme-Mal:93}. For a detailed description of the coordinatization theory for other generalized
polygons we refer to \cite{Mal:98}. The second coordinatization follows in fact from the natural embedding of
the Ree hexagon in $\PG (6,\K)$.
\subsection{Hexagonal sexternary rings for Ree hexagons}
In \cite{Mal:98} a coordinatization theory with respect to a flag $\{(\infty),[\infty]\}$ is described. It is a
generalization of Halls coordinatization for generalized triangles. Here we describe explicitly the
coordinatization of the Ree hexagon. Let $\K$ be a field of characteristic 3 and let $\theta$ be a Tits
endomorphism, i.e., $\theta$ is an endomorphism of $\K$ with the property that $(x^\theta)^\theta =x^3$. Let
$\K'$ be the image of $\K$ under $\theta$ (note that $\K^3 \leq \K' \leq \K$). We define a hexagonal sexternary
ring $\mathcal{R}=(\K,\K',\Psi_1,\Psi_2,\Psi_3,\Psi_4)$ with
$$
\begin{cases}
\Psi_1(k,a,l,a',l',a'')=a^3k+l,\\
\Psi_2(k,a,l,a',l',a'')=a^2k+a'+aa'',\\
\Psi_3(k,a,l,a',l',a'')=a^3k^2+l'+kl,\\
\Psi_4(k,a,l,a',l',a'')=-ak+a'',
\end{cases}
$$
where $a,a',a''\in\K$ and $k,l,l'\in\K'$. This defines the Ree hexagon $\sf{H}(\K,\K')$ as follows.
The points and lines are the $i$-tuples of elements of $\K \cup \K'$ ($i \leq 5$) with alternately an entry in
$\K$ and one in $\K'$, and for points (lines) the last entry is supposed to be in $\K$ ($\K'$), except when
$i=0$, in case we denote the point by $(\infty)$ and the line by $[\infty]$ (we generally use round parentheses
for points and square brackets for lines). Incidence is defined as follows~:
\begin{itemize}
 \item If the number of coordinates of a point $p$ differs by at least 2 from the number of coordinates of a
line $L$, then $p$ and $L$ are not incident. \item If the number $i_p$ of coordinates of a point $p$ differs by
exactly 1 from the number $i_L$ of coordinates of a line $L$, then $p$ is incident with $L$ if and only if $p$
and $L$ share the first $i$ coordinates, where $i$ is the smallest among $i_p$ and $i_L$.
 \item If $i_p = i_L\neq 5$, then $p$ is incident with $L$ if and only if $p=(\infty)$ and $L=[\infty]$.
 \item A point $p$ with coordinates $(a,l,a',l',a'')$ is incident with a line $[k,b,k',b',k'']$
(with above notation, and with $b,b'\in\K$ and $k',k''\in\K'$) if and only if
$$
\begin{cases}
\Psi_1(k,a,l,a',l',a'')=k'',\\
\Psi_2(k,a,l,a',l',a'')=b'',\\
\Psi_3(k,a,l,a',l',a'')=k',\\
\Psi_4(k,a,l,a',l',a'')=b'.
\end{cases}
$$
\end{itemize}

In general, a \emph{generalized hexagon} is a point-line structure with the property that each pair of elements
(en element can be a point, a line or a \emph{flag} | the latter is an incident point-line pair) is contained
in an ordinary hexagon, but never in an ordinary pentagon, quadrangle, triangle or digon. For a general
introduction to generalized hexagons, we refer to \cite{Mal:98}. We mention that generalized hexagons | and
more general, generalized polygons | were introduced by Jacques Tits \cite{Tit:59}.

Roughly, a Ree group is the centralizer of a polarity in a Ree hexagon, and it acts on a Ree-Tits ovoid (namely, the set
of absolute points of the polarity). A \emph{polarity} is an involutive incidence preserving bijection between
points and lines, and between lines and points. A point or line incident with its image is called
\emph{absolute}. An \emph{ovoid} of a generalized hexagon is a set of mutually opposite points (i.e., points at
distance 6 in the incidence graph) such that each point not in the ovoid is collinear to exactly one point of
the ovoid. The dual notion is a \emph{spread}.  It is a general fact that the set of absolute points (lines) of
a polarity is an ovoid (spread), see Chapter 7 of \cite{Mal:98}.

\subsection{The embedding of the Ree hexagon in $\PG (6,\K)$}
The Ree hexagon has a natural embedding in $\PG (6,\K)$.  Indeed, $\sfH(\K,\K')$ is a substructure of the split
Cayley hexagon $\mathsf{H}(\K)$, which has itself a natural embedding in $\PG(6,\K)$ as discovered and
described by Tits in \cite{Tit:59}, see also Chapter 2 of \cite{Mal:98}. Here, we content ourselves with the
table translating the above coordinates to the projective coordinates. We refer to Chapter 3 of \cite{Mal:98}
for details and proofs.

We write $\alpha$ for $-al' +a'^2+a''l+aa'a''$ and $\beta$ for $l-aa'-a^2 a''$.

\bigskip

\begin{tabular}{|l|r|}
\hline
Coordinates in $\sf{H}(\K,\K')$ & Coordinates in $\PG (6,\K)$\\
\hline
$(\infty)$ & $(1,0,0,0,0,0,0)$\\
$(a)$ & $(a,0,0,0,0,0,1)$\\
$(k,b)$&$(b,0,0,0,0,1,-k)$\\
$(a,l,a')$&$(-l-aa',1,0,-a, 0,a^2,-a')$\\
$(k,b,k',b')$& $(k'+bb',k,1,b ,0,b', b^2-b'k)$\\
$(a,l,a',l',a'')$& {\small{$(\alpha,-a'',-a,-a'+aa'',  1, \beta,-l'+a'a'' )$}}\\
\hline
\end{tabular}\\

\begin{tabular}{|l|r|}
\hline
Coordinates in $\sf{H}(\K,\K')$ & Points generating this line\\
\hline
$[\infty]$ & $(\infty)$ and $(0)$\\
$[k]$ & $(\infty)$ and $(k,0)$\\
$[a,l]$&$(a)$ and $(a,l,0)$\\
$[k,b,k']$&$(k,b)$ and $(k,b,0)$\\
$[a,l,a',l']$& $(a,l,a')$ and $(a,l,a',l',0)$\\
$[k,b,k',b',k'']$& $(k,b,k',b')$ and $(0,k'',b',k'+kk'',b)$\\
\hline
\end{tabular}

The subgroup of $\mathsf{PSL}_7(\K)$ stabilizing the point set and line set of $\mathsf{H}(\K,\K')$ is denoted by $\mathsf{G}_2(\K,\K')$ and is simple (a mixed group of type $\mathsf{G}_2$, see \cite{Tit:74}).

\section{The Ree-Tits ovoid}
We start from the Ree hexagon $\sf{H}(\K,\K')$, where $\K'=\K^\theta$, with $\theta$ as above a
Tits-endomorphism of $\K$. This hexagon allows a polarity. The absolute points under this polarity form an
ovoid of the Ree hexagon~: the \emph{Ree-Tits ovoid}, see Chapter 7 of \cite{Mal:98}. We denote the polarity,
which we can choose in such a way that it fixes the flags $\{(\infty),[\infty]\}$ and $\{(0,0,0,0,0),[0,0,0,0,0]\}$
and maps the point $(1)$ onto the line $[1]$, by $\rho$. It has the following actions.
\begin{align*}
&{(a,l,a',l',a'')}^\rho = [a^\theta,l^{\theta^{-1}},a'^\theta,l'^{\theta^{-1}},a''^\theta];\\
&{[k,b,k',b',k'']}^\rho = (k^{\theta^{-1}},b^\theta,k'^{\theta^{-1}},b'^\theta,k''^{\theta^{-1}}).
\end{align*}
for all $a,a',a'',b,b' \in \K$ and $k,k',k'',l,l' \in \K'$.

Now the point $(a,l,a',l',a'')$ is absolute for $\rho$ if and only if it is incident with its image. This leads
to the following conditions~:
$$
\begin{cases}
l= a''^{\theta}-a^{\theta + 3},\\
l'=a^{2\theta + 3}+a'^\theta + a^\theta a''^\theta.
\end{cases}
$$
\textbf{Coordinates of the Ree-Tits ovoid in $\PG(6,\K)$.} | We associate the triple $(a,a'',a'-aa'')$ with the
point $(a,a''^\theta -a^{3+\theta},a',a^{3+2\theta}+a'^\theta +a^\theta a''^\theta,a'' )$. Now, for $a,a',a''
\in \K$, we put
\begin{align*}
&f_1(a,a',a'')= -a^{4+ 2\theta}-aa''^\theta +a^{1+\theta}a'^\theta +a''^2 +a'^{1+\theta} -a'a^{3+\theta}-a^2a'^2,\\
&f_2(a,a',a'')=-a^{3+\theta}+a'^\theta -aa''+a^2a',\\
&f_3(a,a',a'')=-a^{3+2\theta}-a''^\theta +a^\theta a'^\theta +a'a''+aa'^2.
\end{align*}
So the set of absolute points can be described in $\PG(6,\K)$ by
{\small{
\begin{align*}
\Omega=&\{(1,0,0,0,0,0,0)\}\cup \\&\{(f_1(a,a',a''),-a',-a,-a'',1,f_2(a,a',a''),f_3(a,a',a''))~\vert~ a,a',a''
\in \K\}.
\end{align*}}}

\bigskip

\textbf{Compact notation.} | As before, we associate the triple $(a,a'',a'-aa'')$ with the point $(a,a''^\theta
-a^{3+\theta},a',a^{3+2\theta}+a'^\theta +a^\theta a''^\theta )$. The set of absolute points under the polarity
is now
$$\Omega=\{(\infty)\} \cup \{(a,a',a'')~|~a,a',a''\in\K\}.$$
The elements of the \emph{root group} $U_\infty$ (a root group is a unipotent subgroup) of the ovoid fixing the
point $(\infty)$ act as follows on the remaining points $(x,x',x'')$~: the unipotent element that fixes
$(\infty)$ and maps $(0,0,0)$ to $(y,y',y'')$ maps $(x,x',x'')$ to $$(x,x',x'') \cdot
(y,y',y'')=(x+y,x'+y'+xy^\theta,x''+y''+xy'-x'y-xy^{\theta + 1}),$$ and this action can also be seen as the
multiplication inside $U_\infty$, see Chapter 7 of \cite{Mal:98}.

We obtain the Ree-Tits Moufang set. Formally, a \emph{Moufang set} is a set $X$ together with a family of
(root) groups $(U_x)_{x\in X}$ such that $U_x$ acts sharply transitively on $X\setminus\{x\}$ and fixes $x$,
and such that each $U_y$, $y\in X$, acts by conjugation on the set $\{U_x~|~x\in X\}$. The (simple) Ree groups
arise as (simple subgroups of the) centralizers of polarities in these hexagons. More exactly, the Ree group $\mathsf{R}(\K,\theta)$ is defined as the centralizer in $\mathsf{G}_2(\K,\K')$ of the outer automorphism $\rho$. This group is simple if $|\K|>3$ and the multiplicative group of $\K$ is generated by all squares together with $-1$, see \cite{Ree:61}. In any case, the group generated by the root groups is simple, provided $|\K|>3$, and it coincides with the derived group $\mathsf{R}'(\K,\theta)$. For $|\K|=3$, $\mathsf{R}(\K,\theta)=\mathsf{R}(3)$ is isomorphic to $\mathsf{P\Gamma L}_2(8)$ and contains $\mathsf{PSL}_2(8)$ as a simple subgroup of index $3$.

We can see the Ree-Tits ovoid
and its automorphism group embedded in the Ree hexagon as a representation of the Ree-Tits Moufang set.
Henceforth, we will denote by $\Omega$ the Ree-Tits ovoid, and by $U_x$, $x\in \Omega$, the root group fixing
$x$ in the Ree-Tits Moufang set over the field $\K$ with associated Tits endomorphism $\theta$.

We will also need the explicit form of a generic element of the root group $U_{(0,0,0)}$, which we shall
briefly denote by $U_0$. This is best given by the action on coordinates in the projective space. Such generic
element $u^{(0,0,0)}_{(x,x',x'')}$ then looks like (and $x,x',x''$ are arbitrary in $\K$):

$\bar{x}=(x_0~x_1~x_2~x_3~x_4~x_5~x_6)\mapsto$
{\small{$$\bar{x}\cdot\left(\begin{array}{ccccccc}
1&f_2(x,x',x'')&f_3(x,x',x'')& x'' &             f_1(x,x',x'')&-x'&-x \\
0 & 1 & -x^\theta            & 0 &             x'-x^{1+\theta} & 0 & 0  \\
0 & 0 & 1                    & 0 &  x & 0 & 0 \\
0 & -x & x'                 & 1 &              -x'' & 0 & 0\\
0&0&0                        & 0&             1 & 0 & 0 \\
0 & x^2 & -x''-xx'           & x &              p & 1 & 0 \\
0 & r & s                    & -x'+x^{1+\theta} &              q & x^\theta & 1
\end{array}\right),$$}}

where $$\begin{cases} p=x^{3+\theta}-{x'}^\theta-xx''-x^2x',\\ q={x''}^\theta + x^\theta{x'}^\theta - x{x'}^2 -
x^{2+\theta}x' - x^{1+\theta}x'' - x^{3+2\theta},\\ r=x''-xx'+x^{2+\theta},\\ s = {x'}^2 -
x^{1+\theta}x'-x^\theta x'',\end{cases}$$

see Section 9.2.4 of \cite{Tha-Tha-Mal:06}.

We are now ready to define the Ree geometries.

\section{The Ree geometry}\label{Reegeometry}
As already mentioned, the Ree groups have root groups of nilpotency class 3. As a consequence, the geometries
that we will define corresponding to the Ree groups will have rank 3. This means that we will have two types of
\emph{blocks} in our geometry and that blocks of one type are subsets of the others. In order to distinguish
the two types of blocks, we will call the ``smallest'' ones \emph{circles}, and the others \emph{spheres}.

The point set $\cP$ of our Ree geometry $\mathcal{G}=(\cP,\cB)$ is the Ree-Tits ovoid. The circles arise as
orbits of a point $y$ under the center $Z(U_x)$ for some point $x \in \cP\setminus\{y\}$, together with that
point $x$. This particular point $x$ is then called a \emph{gnarl} of this circle. So every point and gnarl
defines a circle in a unique way. The spheres are again a point $x$ together with the orbit of some point $y$,
$y\neq x$, but this time under the group $U_x'=[U_x,U_x]$. The point $x$ is a \emph{gnarl} of the sphere. The
circles and spheres together form the block set $\cB$ of $\mathcal{G}$. The set of circles will be denoted by
$\cC$, and the set of spheres by $\cS$. We can define two further geometries by restricting the set of blocks.
We call the geometries $\mathcal{G_C}=(\cP,\cC)$ and $\mathcal{G_S}=(\cP,\cS)$ the \emph{truncated Ree
geometries}. The gnarls of circles and spheres will turn out to be unique.

Let us be more concrete now and look for the coordinates of the circles and spheres which have $(\infty)$ for
gnarl.


We first claim that, if $|\K|>3$, then the group $U'_\infty$ is precisely $\{(0,u',u'')~\vert~ u',u'' \in
\K\}$. Indeed, computing an arbitrary commutator, we get \begin{align*}
[(u_1,u'_1,u''_1),(u_2,u'_2,u''_2)]&=(0,u_1u_2^\theta-u_2u_1^\theta,u_1'u_2-u_1u'_2-u_1u_2^{1+\theta}+u_2u_1^{1+\theta}).
\end{align*}

Noting that $(0,x',0)\cdot(0,0,x'')=(0,x',x'')$, we only have to show that $(0,x',0)\in U'_\infty$, for all
$x'\in\K$, and that $(0,0,x'')\in U'_\infty$, for all $x''\in\K$. Putting $u_1=u_1''=u_2'=u_2''=0$, $u_1'=1$
and $u_2=x''$ in the above commutator, we see that $(0,0,x'')\in U_\infty$. Now let $x'\in \K$ be arbitrary.
Since $|\K|>3$, there exists an element $t\in\K$ with $t^3-t\neq 0$. Put $k=t^3-t$ and let $y=x'k^{-\theta}$.
Putting $u_1'=u_2'=u_1''=u_2''=0$ and $(u_1,u_2)=(y,t^3)$, respectively $(u_1,u_2)=(t^\theta y,1)$, we obtain
$(0,t^{3\theta}y-t^3y^\theta,0)\in U_{\infty}'$ and $(0,t^\theta y-t^3 y^\theta,0)\in U_\infty'$. Multiplying
the former with the inverse of the latter, we see that $(0,x',0)\in U_\infty'$, proving our claim.

If $|\K|=3$, then $U_\infty'$ has order 3 and coincides with the center (see below). In this case, we will
substitute $U'_x$ by the subgroup of $U_x$ generated by the elements of order $3$, and we will denote it, with
abuse of notation, by $U'_x$ (but there will be no confusion possible), since for $|\K|>3$, the derived group
coincides with the group generated by elements of order $3$ (as one can check easily).

The center of $U_\infty$ is the subgroup $\{(0,0,u'')~\vert~ u'' \in \K \}$. Indeed, this follows from the
explicit form of the multiplication in $U_\infty$ by standard arguments.  Since the commutator of an element
$(0,u_1',u_1'') \in U'_\infty$ and $(u_2,u_2',u_2'') \in U_\infty$ is
\begin{align*}
[(0,u'_1,u''_1),(u_2,u'_2,u''_2)]&=(0,0,u_1'u_2)\\
&=(0,0,u''),
\end{align*}
with $u''$ essentially arbitrary, we see that the second derived group $U_\infty''$ coincides with the center
$Z(U_\infty)$.

Now, since the circles having $(\infty)$ as gnarl are the orbits of a point $(a,a',a'')$ under the group
$\{(0,0,x)~\vert~ x \in \K\}$, all circles with $(\infty)$ as gnarl are given by
$$\{(a,a',a''+x)~\vert~ x \in \K\}\cup \{(\infty)\}=\{(a,a',t)~\vert~ t \in \K\}\cup \{(\infty)\}.$$
The spheres with gnarl $(\infty)$ have the following description:
$$\{a,a'+x',a''+x''+ax'~\vert~ x',x'' \in \K \} \cup\{(\infty)\}=\{(a,t',t'')~\vert~ t', t'' \in  \K\} \cup \{(\infty)\}.$$

We can now interpret the algebraic description of a circle and a sphere with gnarl $(\infty)$
in the corresponding Ree hexagon $\sf{H}(\K,\K')$. We leave it as an easy exercise to the
reader to see that these circles and spheres look as follows~: a circle is the set of absolute points at
distance 3 from a non-absolute line $M$, not going through an absolute point. The unique absolute point for
which its corresponding absolute line intersects $M$ is the gnarl of the circle. Likewise, a sphere is the set
of absolute points not opposite some point $p$, with $p$ lying on an absolute line. The unique absolute point
at distance 2 from $p$ is the gnarl of the sphere. Conversely, every such set is a circle or sphere,
respectively. It follows now easily that the gnarl of a circle and of a sphere is unique. These gnarls will
play a prominent role in our proofs.


As an application we make the following important observation.

\begin{lemma}\label{samegnarl}
A sphere contains only circles with the same gnarl. Also, the point set of a sphere, except for its gnarl, is
partitioned by the circles contained in the sphere.
\end{lemma}
\proof Let us consider a sphere and circle, and assume that this sphere's gnarl is the absolute point $p$ while
the gnarl of the circle is a different absolute point $q$. The flags $\{p,p^\rho\}$ and $\{q,q^\rho\}$ determine an
unique apartment $\Sigma$ containing both flags, and because both flags are absolute $\rho$ will stabilize
$\Sigma$. Denote the unique line in $\Sigma$ at distance 2 from $q^\rho$ and at distance 3 from $p$ with $L$
and the projection of $q$ on $p^\rho$ with $r$. Let $a$ be a third absolute point on the circle different from
both $p$ and $q$.

Because $a$ lies on the circle with gnarl $q$ through $p$, $a$ lies at distance 3 from $L$. Similarly $a$ also
lies at distance 4 from $r$ because of the definition of a sphere. The last statement implies that $a^\rho$
lies at distance 4 from the line $r^\rho$. This line $r^\rho$ intersects the line $L$, so the point $a$ and the
lines $L, r^\rho, a^\rho$ are contained in a ordinary 5-gon, which contradicts the definition of a generalized
hexagon. This proves the first assertion.

For the second assertion, we just consider the circles defined by the non-absolute lines of
$\sfH(\K,\K^\theta)$ through the point defining the sphere in question. \qed

\section{Statement of main results}\label{mainresult}
Given the construction of the circles and spheres in the corresponding Ree hexagon $\sfH(\K,\K^\theta)$, it is
clear that every collineation of $\sfH(\K,\K^\theta)$ that commutes with the polarity $\rho$ induces a
collineation of the Ree geometry and its truncations. Our main results now say that also the converse holds.
More precisely~:

\textbf{Main Result 1.} \emph{The full automorphism group of the Ree geometry $\cG=(\cP,\cB)$ is the
centralizer of $\rho$ in the full collineation group of $\sfH(\K,\K^\theta)$.}

Likewise, we will show:

\textbf{Main Result 2.} \emph{The full automorphism groups of the truncated Ree geometries
$\mathcal{G_C}=(\cP,\cC)$ and $\mathcal{G_S}=(\cP,\cS)$ coincide with the centralizer of $\rho$ in the full
collineation group of $\sfH(\K,\K^\theta)$.}

As a main consequence we will be able to show:

\textbf{Main Corollary.} \emph{The stabilizer of a Ree-Tits ovoid in the full collineation group of
$\sfH(\K,\K^\theta)$ coincides with the centralizer of the corresponding polarity in the full collineation
group of $\sfH(\K,\K^\theta)$. Consequently, any polarity is determined by its set of absolute points.}

The latter was already announced in \cite{Mal:98} as Theorem~7.7.9, but not proved there.

We will now prove these results.

\section{Auxiliary tools}

\subsection{The derived geometry at $(\infty)$}\label{derived}

We define the structure $\mathcal{G'}=(\cP',\cB')$, where $\cP'=\cP \setminus \{(\infty)\}$, and $\cB'$ is the
set of blocks of $\mathcal{G}$ going through $(\infty)$, with $(\infty)$ removed. We call this the
\emph{derived geometry at $(\infty)$}, inspired by a similar concept in the theory of designs. In order to know
the coordinates of the circles through $(\infty)$ we first write down the coordinates of the circles with gnarl
$(\infty)$. As we saw earlier, these are the sets
\[\{(a,a',t)~\vert~ t \in \K\}\cup (\infty),\mbox{ with }a,a'\in\K.
\]
Removing the point $(\infty)$ gives us the \emph{vertical line} $L_{a,a'}$. We now compute the coordinates of
the circle with gnarl $(0,0,0)$ through $(\infty)$. The point $(\infty)$ is identified with $(1,0,0,0,0,0,0)$,
so its orbit under $Z(U_0)$ (using the elements $u^{(0,0,0)}_{(0,0,x'')}$ defined above) is the set

\begin{align*} \{(1,f_2(0,0,x''),f_3(0,0,x''),& x'',f_1(0,0,x''),0,0)~|~x''\in\K\}\\&=\{(1,0,-x''^\theta
,x'',x''^2,0,0)~|~x''\in\K\}. \end{align*}

Putting $x={x''}^{-2-\theta}$ (and hence $x''=x^{-2+\theta}$), adding the gnarl and deleting the point
$(\infty)$, we obtain the set $\{(x,0,-x^{2+\theta})~|~x\in\K\}$. The image of this set under $(a,a',a'') \in
U_\infty$ is the set
\[
\{(a+x,a'+a^\theta x,a''+(a'-a^{1+\theta})x-x^{2+\theta})~\vert~ x \in \K\},
\]
which we call the \emph{ordinary line} $C_{(a,a',a'')}$ (\emph{with gnarl} $(a,a',a'')$). Note that unlike the vertical lines, these are not affine lines.

Just as we did for circles, we consider the spheres with gnarl $(\infty)$ and the other spheres through
$(\infty)$ separately.

The spheres with gnarl $(\infty)$ are the sets $\{(a,t',t'')~\vert~ t', t'' \in  \K\} \cup \{(\infty)\}$, with
$a\in\K$. Removing the point $(\infty)$ gives us the \emph{vertical plane} $P_a$.

The orbit of $(\infty)$ under $U_0'$, using the elements $u^{(0,0,0)}_{(0,x',x'')}$ is the set {\small{
\begin{align*}
\{(1,f_2(0,x',&x''),f_3(0,x',x''),x'',f_1(0,x',x''),-x',0)~|~x',x''\in\K\}\\&=\{(1,x'^\theta ,-x''^\theta
+x'x'',x'',x''^2+x'^{1+\theta},-x',0)~|~x',x''\in\K\}\\&=\left\{\left(\frac{x''^\theta
-x'x''}{x''^2+x'^{1+\theta}},\frac{-x'^\theta}{x''^2+x'^{1+\theta}},\frac{-x''}{x''^2+x'^{1+\theta}}\right)~|~\K\times\K\ni
(x',x'')\neq(0,0)\right\} \cup \{(\infty)\}.
\end{align*}
}} Note that ${x''}^2\neq-{x'}^{1+\theta}$ is equivalent with $(x',x'')\neq(0,0)$. Adding $(0,0,0)$, the image
of this sphere under $(a,a',a'') \in U_\infty$ is the set {\small{
\begin{align*}
\biggl\{\biggl(\frac{x''^\theta -x'x''}{x''^2+x'^{1+\theta}},\frac{-x'^\theta}{x''^2+x'^{1+\theta}},
\frac{-x''}{x''^2+x'^{1+\theta}}\biggr)\cdot(a,a',a'')~|~\K\times\K\ni(x',x'')&\ne(0,0)\biggr\} \\ &\cup
\{(a,a',a''),(\infty)\}. \end{align*}}}

Removing the point $(\infty)$ gives us the \emph{ordinary plane} $S_{(a,a',a'')}$ (\emph{with gnarl} $(a,a',a'')$). Again note that these are not affine planes, while the vertical planes are.

Notice that points of vertical planes have constant first coordinate, while the points of an ordinary line
never have constant first coordinate. This provides an algebraic proof of Lemma~\ref{samegnarl}.

\subsection{Parallellism in this derived structure}\label{subsection:para}

We consider the set of points $(x,x',x'')$ as an affine space in the standard way, and call the planes
\emph{affine planes}.  We assume that the coordinates are given with respect to a basis with axes $X,Y,Z$.

First we remark that every ordinary line $C_{(a,a',a'')}$ completely lies in the affine plane with equation
$Y=a^\theta X+(a'-a^{1+\theta})$. We say that two ordinary lines $C_1$ and $C_2$ are \emph{parallel} if all
vertical lines intersecting $C_1$ intersect $C_2$ | in that case the two ordinary lines lie in the same affine
plane | or if there is no vertical line intersecting both ordinary lines | which implies that the ordinary
lines lie in parallel, but disjoint, affine planes.


We claim that two ordinary lines $C_{(a,a',a'')}$ and $C_{(b,b',b'')}$ are parallel if and only if $a=b$.
Indeed, a vertical line meeting the ordinary line $C_{(a,a',a'')}$ must lie in the affine plane $Y=a^\theta
X+(a'-a^{1+\theta})$, so any vertical line meeting both $C_{(a,a',a'')}$ and $C_{(b,b',b'')}$ must lie in the
intersection of
\[
\begin{cases}
Y=a^\theta X+(a'-a^{1+\theta}),\\
Y=b^\theta X+(b'-b^{1+\theta}).
\end{cases}\]
This has a unique solution if and only if $a \neq b$, proving our claim.


We have the following direct lemma.

\begin{lemma}\label{cor:sfeeruitcirkels}
The gnarls of the ordinary lines of the parallel class of $C_{(a,a',a'')}$ are
exactly the points of the vertical plane
$P_a$.
\end{lemma}
\proof The above says that the set of gnarls of the lines of the parallel class of $C_{(a,a',a'')}$ is given by
$\{(a,t',t'')~\vert~ t', t'' \in  \K\}$, which is exactly $P_a$. \qed

\subsection{Ree unitals}
In Section~\ref{section:sferen}, we will use the Ree unitals mentioned in the introduction. We do not need a formal definition, or a
complete description of them, but only the following two facts~: (1) two different points are joined by exactly
one block of the Ree unital, (2) the block through $(\infty)$ and $(a,0,a'')$, with $a$ and $a'' \in \K$, is
given by $\{(\infty)\} \cup \{(a,t,a''-at)|t\in \K\}$ (see Chapter 7 of \cite{Mal:98}).

If $B$ is a unital block containing $(\infty)$, then we will call the set $B\setminus\{(\infty)\}$ an
\emph{affine unital block}.

\section{Automorphism group of the Ree geometry}
\textbf{General idea.} | We consider an automorphism $\varphi$ of the Ree geometry. Without loss of generality
we may assume that $\varphi$ fixes both $(\infty)$ and $(0,0,0)$. We will prove that $\varphi$ must preserve
gnarls, and this will imply that it has to preserve the parallelism we just defined. We then compute the
algebraic form of $\varphi$ and conclude that it can be extended to $\sfH(\K,\K^\theta)$.

\begin{lemma}\label{gnarls}
The automorphism $\varphi$ maps the gnarl of any sphere onto the gnarl of the image of the sphere, and it maps
the gnarl of any circle onto the gnarl of the image of the circle under $\varphi$.
\end{lemma}
\proof Any automorphism of $\Delta$ maps spheres onto spheres and circles onto circles, since every circle is
properly contained in a sphere, but no sphere is properly contained in any circle or sphere. Since the gnarl of
a sphere is exactly the intersection of all circles contained in it (by Lemma~\ref{samegnarl}), and there are
at least two such circles, $\varphi$ preserves gnarls of spheres. But then $\varphi$ must also preserve the
gnarls of these circles. \qed

Since $\varphi$ fixes the points $(\infty)$ and $(0,0,0)$, it acts on the derived structure $\cG'$, and the
previous lemma implies that $\varphi$ fixes the set of vertical lines. Therefore the points $(a,a',z_1)$ and
$(a,a',z_2)$ are mapped on the same vertical line. If we represent $\varphi$ as follows~:\[\varphi :
(x,y,z)\mapsto (g_1(x,y,z),g_2(x,y,z),g_3(x,y,z))\] then both $g_1$ and $g_2$ have to be independent of $z$ and
we write $g_i(x,y,z)=g_i(x,y)$, $i=1,2$.

The mapping $\varphi$ preserves the parallel relation between ordinary lines, since the number of vertical
lines meeting two circles (i.e. none, one or all) is preserved under $\varphi$. This translates to $g_1$ being
independent of $y$. Indeed, two points $(a,y_1,z_1)$ and $(a,y_2,z_2)$ being the gnarls of two parallel
ordinary lines are mapped onto two gnarls of parallel ordinary lines, which implies that $
g_1(a,y_1)=g_1(a,y_2)$ for every choice for $y_1$ and $y_2$.

The point $(0,0,0)$ is fixed by $\varphi$, so the affine plane $Y=0$ | which is the unique affine plane
containing both $C_{(0,0,0)}$ and $L_{0,0}$, and which consist of the union of vertical lines all meeting
$C_{(0,0,0)}$ | is fixed by $\varphi$. The plane $Y=c_1$ | which is also a union of vertical lines | must
necessarily get mapped onto a plane $Y=c_2$. So $g_2(x,c_1)=g_2(0,c_1)$ for every choice of $x\in\K$.

It follows that there are two permutations $\alpha$ and $\beta$ of $\K$ such that $(x,y,z)^\varphi$ is equal to
$(x^\alpha,y^\beta,g_3(x,y,z))$. Since $\varphi$ preserves gnarls, it maps the ordinary line $C_{(a,b,c)}$ onto
the ordinary line $C_{(a^\alpha,b^\beta,g_3(a,b,c))}$. Now notice that the point $(x,y,z)$ can only be
contained in the ordinary line $C_{(a,b,c)}$ if $y=b+a^\theta (x-a)$. Expressing that the point $(a+x,y,z)$
lies on the circle $C_{(a,b,c)}$ if and only if its image under $\varphi$ lies in $C_{(a,b,c)}^{\varphi}$ shows
that, for all $a,b,x\in \K$,
\begin{equation}(b+a^\theta x)^\beta=b^\beta+(a^\alpha)^\theta ((x+a)^\alpha-a^\alpha).\label{g2}\end{equation}
Putting $b=0$, and noting that $0^\alpha=0^\beta=0$, we see that $(a^\alpha)^\theta
((x+a)^\alpha-a^\alpha)=(a^\theta x)^\beta$, which implies, by substituting this back in Equation~(\ref{g2}),
that $(b+a^\theta x)^\beta=b^\beta+(a^\theta x)^\beta$. So $\beta$ is additive.  Put $\ell=1^\alpha$. Then we
see, by setting $a=1$ and $b=0$ in the Equation~(\ref{g2}) above, that
\begin{equation}x^\beta=\ell^\theta ((x+1)^\alpha-1^\alpha), \label{beta}\end{equation}
so $\alpha$ is additive if and only if $(x+1)^\alpha=x^\alpha+1^\alpha$. Plugging in $x=m-1$ in Equation~(\ref
{beta}) we have that $(m-1)^\beta =\ell^\theta (m^\alpha -1^\alpha)$. Because of the additivity of $\beta$ we
have on the other hand that $(m-1)^\beta = m^\beta + (-1)^\beta=\ell^\theta ((1+m)^\alpha -2 \cdot 1^\alpha)$.
So $\alpha$ is additive as well.

We now have that $x^\beta=\ell^\theta x^\alpha$.
We can define the bijection $\sigma:\K\rightarrow \ell^{-1}\K:y\mapsto y^\sigma=\ell^{-1}y^\alpha$
(note that $1^\sigma=1$). Plugging in these identities in equation (\ref{g2}) yields
$$(b+a^\theta x)^\sigma=b^\sigma+(a^\sigma)^\theta x^\sigma,$$
for all $a,b,x\in \K$. Putting $a=1$, we see that $\sigma$ is additive; putting $b=0$ and $x=1$, we see that
$\sigma$ commutes with $\theta$. Putting $b=0$, we see that $(xy)^\sigma=x^\sigma y^\sigma$ for $x\in
\K^\theta$ and $y\in \K$. If $x,y\in\K$, then $$((xy)^\sigma)^\theta=((xy)^\theta)^\sigma=(x^\theta
y^\theta)^\sigma=(x^\theta)^\sigma (y^\theta)^\sigma=(x^\sigma)^\theta(y^\sigma)^\theta=(x^\sigma
y^\sigma)^\theta,$$ and the injectivity of $\theta$ implies that $\sigma$ is an automorphism of $\K$. Now the
action of $\varphi$ on a point $(x,y,z)$ is given by $(x,y,z)^\varphi=(\ell x^\sigma,\ell^{1+\theta}
y^\sigma,g_3(x,y,z))$, for all $x,y,z\in \K$.

Let us now investigate what $g_3(x,y,z)$ looks like.

The point $p$ with coordinates {\small
$(a-\frac{a'}{a^\theta},0,a''+(a'-a^{1+\theta})(\frac{-a'}{a^\theta})-(\frac{-a'}{a^\theta})^{2+\theta})$} lies
on both $C_{(a,a',a'')}$ and on the ordinary line with gnarl
$(0,0,a''+\frac{(a^{1+\theta}-a')^{1+\theta}+a'^{1+\theta}}{a^{2+\theta}})$. So its image under $\varphi$ lies
on the ordinary line with gnarl $(\ell a^\sigma,\ell^{1+\theta}a'^\sigma,g_3(a,a',a''))$ and on the ordinary
line with gnarl $(0,0,g_3(0,0,a''+\frac{(a^{1+\theta}-a')^{1+\theta}+a'^{1+\theta}}{a^{2+\theta}}))$. This
leads to {\footnotesize
\[\begin{cases}
g_3(a-\frac{a'}{a^\theta},0,a''-\frac{(a'-a^{1+\theta})a'}{a^\theta}+(\frac{a'}{a^\theta})^{2+\theta})=
g_3(a,a',a'')-\ell^{2+\theta}(\frac{a'^2}{a^\theta}-aa'-\frac{a'^{2+\theta}}{a^{3+2\theta}})^\sigma, \\
\\
g_3(a-\frac{a'}{a^\theta},0,a''-\frac{(a'-a^{1+\theta})a'}{a^\theta}+(\frac{a'}{a^\theta})^{2+\theta})=
g_3(0,0,a''+\frac{(a^{1+\theta}-a')^{1+\theta}+a'^{1+\theta}}{a^{2+\theta}})-
(\ell(a-\frac{a'}{a^\theta})^\sigma)^{2+\theta}.
\end{cases}
\]}

Putting these two equations together we get~:

\[
g_3(a,a',a'')=g_3\biggl(0,0,a''+\frac{(a'-a^{1+\theta})^{1+\theta}+a'^{1+\theta}}{a^{2+\theta}}\biggr)-
\ell^{2+\theta}\biggl(\frac{(a'-a^{1+\theta})^{1+\theta}+a'^{1+\theta}}{a^{2+\theta}}\biggr)^\sigma,
\]
for every $a \in \K \backslash \{0\}$ and $a',a'' \in \K$. We want to extend this equation to one with $a=0$.
To this end, we note that the point $(0,a',a'')$ lies on every circle with gnarl
$(A,a'+A^{1+\theta},a''+a'A-A^{2+\theta})$, with $A\in\K$. We now only consider $A\neq 0$. Then we take the
image under $\varphi$ and obtain that $$g_3(0,a',a'')=g_3(A,a'+A^{1+\theta},
a''+a'A-A^{2+\theta})-\ell^{2+\theta}(Aa'-A^{2+\theta})^\sigma.$$ We can now use the above expression for
$g_3(a,a',a'')$ for $a\neq 0$ to express $g_3(0,a',a'')$ in terms of $g_3(0,0,z)$, for some $z\in\K$. We carry
this out explicitly, and substitute $a'=B^{\theta-1}$ and $A=B^{2-\theta}$, and obtain after a tedious
calculation

$$g_3(0,B^{\theta-1},a'')=g_3(0,0,a''-B)+\ell^{2+\theta}B^\sigma,$$

for all $B\in\K\setminus\{0\}$, and all $a''\in\K$. Substituting $-B$ for $B$, we see that
$g_3(0,0,a''-B)=g_3(0,0,a''+B)+\ell^{2+\theta}B^\sigma$. We may now put $a''=-B$ and obtain finally that
$g_3(0,0,B)=\ell^{2+\theta}B^\sigma$.  Plugging this into the formulae above for $g_3(a,a',a'')$, $a\neq 0$,
and $g_2(0,a',a'')$, we see that $g_3(a,a',a'')=\ell^{2+\theta}{a''}^\sigma$, for all $a,a',a''\in\K$.

So the action of $\varphi$ on a point $(x,y,z)$ is given by $(x,y,z)^\varphi=(\ell
x^\sigma,\ell^{1+\theta}y^\sigma,\ell^{2+\theta} z^\sigma)$, with $\sigma$ and $\theta$ commuting automorphisms
of $\K$. This action is the restriction to $\Omega$ of the collineation of $\sfH(\K,\K^\theta)$ defined by the
following mapping on the points and lines with five coordinates:

$$\begin{cases}
(a,l,a',l',a'')\mapsto (\ell a^\sigma, \ell^{\theta+3}l^\sigma,\ell^{\theta+2}{a'}^\sigma,
\ell^{2\theta+3}{l'}^\sigma,\ell^{\theta+1}{a''}^\sigma),\\
[k,b,k',b',k'']\mapsto [\ell^\theta k^\sigma, \ell^{\theta+1}b^\sigma,\ell^{2\theta+3}{k'}^\sigma,
\ell^{\theta+2}{b'}^\sigma,\ell^{\theta+3}{k''}^\sigma].
\end{cases}
$$

The proof is complete. \qed

\section{Automorphism group of the truncated Ree geometry restricted to points and circles}
\textbf{General idea.} | Let $\mathcal{G_C}=(\cP,\cC)$ be the truncated Ree geometry, with $\cC$ the set of
circles. We first prove that the gnarls of circles have to be mapped onto gnarls of circles. Then we use the
result from the previous section to prove that the automorphism group of $\mathcal{G_C}$ is equal to the
automorphism group of the Ree geometry $\mathcal{G}$.

We denote by $\mathcal{G_C}'$ the derived geometry in $(\infty)$ (so the point set is
$\cP\setminus\{(\infty)\}$ and the blocks are the vertical and ordinary lines, as defined in
Section~\ref{derived}).

\begin{lemma}\label{gnarlscircles}
The full group $G$ of automorphisms of $\mathcal{G_C}'$ has two orbits on the lines, which are the vertical and
the ordinary lines.
\end{lemma}
\proof It is clear that $G$ acts transitively on both the set of vertical lines and the set of ordinary lines
(as $G$ contains the corresponding Ree group), so we only have to exclude the possibility of one orbit. We
suppose this is the case and derive a contradiction.

Consider, as previously, the point set $\cP\setminus\{(\infty)\}$ as a 3-dimensional affine space with point
set $\{(a,a',a'')|a,a',a'' \in \K \}$. We project it on the 2-dimensional space $\{(a,a',0)|a,a' \in \K \}$ by
the standard projection map $(a,a',a'') \mapsto (a,a',0)$. The projection of a vertical line $L_{a,a'}$ is the
point $(a,a',0)$, and the projection of an ordinary line $C_{(a,a',a'')}$ is the affine line $Y=a^\theta X +
(a' - a^{1+\theta})$. The images of a parallel class of ordinary lines will form a parallel class of affine lines in the 2-dimensional space. This implies that all these affine lines form the line set of a net $\cN$, and a parallel class of ordinary
lines is projected to a parallel class in this net. (A net is a point-line geometry where for each point $p$ not on a line $L$, there exist an unique line incident with $p$, parallel with $L$).

Let $L$ be a vertical line and $M$ a vertical or ordinary line disjoint from $L$. If $M$ is a vertical line
then the projection of $L$ and $M$ are two points. If there exists one ordinary line such that the projection
contains both points, then translating this back to the lines means through each point of $L$ there is an
(ordinary) line intersecting $M$ (by varying the third coordinate $a''$). If, on the other hand, there is no
projection of an ordinary line containing both points, then there is no (ordinary) line intersecting both $L$
and $M$.

If $M$ is an ordinary line, then the projection of $M$ is a certain affine line with equation $Y=a^\theta X +
(a' - a^{1+\theta})$. As no projection of an ordinary line is of the form $X = c$ with $c\in \K$ a constant,
there are points of $M$ through which no (ordinary) line passes that also intersects $L$ (because we would have
projections of the form $X = c$). Also, there obviously are ordinary lines whose projection contains the
projection of $L$ and intersect the projection of $M$. The set of ordinary lines projected to this projection
forms a subset of a parallel class exactly one member of which intersects both $L$ and $M$. We conclude that
there exist lines intersecting both $L$ and $M$, but not through each point of $M$.

In the above two paragraphs we proved that we can tell a vertical line from an ordinary line if one vertical
line is given. Using the hypothesis that there is only one orbit on the lines, this implies that there is an
equivalence relation on the lines which is preserved by $G$.  One of the equivalence classes is obviously the
set of vertical lines. By transitivity it follows that through each point of $\mathcal{G_C}'$ there is exactly
one line of a given equivalence class. We now claim that the other classes are the parallel classes of ordinary
lines. Indeed, if an ordinary line $C_{(a,a',a'')}$ lies in a certain equivalence class, then all lines
$C_{(a,a',k)}$ with $k \in \K$ lie in this class, because there is a vertical line through each point of
$C_{(a,a',a'')}$ intersecting $C_{(a,a',k)}$. The subsets of the equivalence classes thus obtained are exactly
the set of lines projected to a common affine line. As two intersecting affine lines can be viewed as the
projection of two intersecting ordinary lines, two of these subsets are parallel if and only if the
corresponding affine lines are parallel. This implies that the equivalence classes are subpartitions of the
parallel classes. But since through each point there has to be a line of each equivalence class, the latter
must coincide with a parallel class.

Now consider the ordinary line $C_{(0,0,0)}$ and its parallel class $\pi$. We can conjugate the center of
$U_{(\infty)}$ to obtain  an automorphism $\phi\in G$ that fixes the ordinary lines in $\pi$, acts freely on
the points of such a line, fixes the equivalency classes, and maps $(0,0,0)$ to $(1,0,-1)$.

Let $(x,x',x'')$ be an arbitrary point of $\mathcal{G_C}'$. This point lies on the ordinary line $C_{(0,x',b)}
= \{(t,x',b + x't-t^{2+\theta})~|~ t \in \K \}$ for $t=x$ with $b := x''-x'x+x^{2+\theta}$. As this ordinary
line is an element of $\pi$, the point $(x,x',x'')^\phi$ also lies on this line. Hence there exists an
$f_{x',b}(x)\in\K$ such that $(x,x',x'')^\phi = (f_{x',b}(x),x',b + x'f_{x',b}(x)-f_{x',b}(x)^{2+\theta})$.
Notice that the middle coordinate is always fixed.

The vertical line $L_{x,x'} = \{(x,x',t) ~|~ t \in \K \}$ must be mapped to another vertical line
$L_{f_{x',b}(x),x'} = \{(f_{x',b}(x),x',t)~|~t \in \K\}$. From this it follows that the function $f$ is
independent from the the last coordinate. As both the first and second coordinate are independent from the
last, it follows that $\phi$ induces an automorphism $\phi'$ on the net $\cN$, mapping $(x,x',0)$ to
$(f_{x',b}(x), x',0)$. Now $\phi'$ also fixes every parallel class of $\cN$ (the parallel class coming from
$\pi$ is even fixed linewise) and maps $(0,0,0)$ to $(1,0,0)$ (because $(0,0,0)^\phi = (1,0,-1)$). It is now
easy to see that this implies $f_{x',b}(x) = x+1$. This gives us the following explicit formula for $\phi$~:
\begin{align}
\phi : (x,x',x'') \mapsto &(x+1,x',x'' - x'x+x^{2+\theta} +x'(x+1) - (x+1)^{2+\theta}) \\
\mapsto &(x+1,x',x'' +x' +x^{2+\theta} - (x+1)^{2+\theta}).
\end{align}
The image of the ordinary line $C_{(1,1,0)} = \{(1+t,1+t,-t^{2+\theta}) ~|~ t \in \K \}$, using the formula for
$\phi$, is~:
\begin{equation}
C_{(1,1,0)}^\phi = \{(t-1,t+1,-t^{2+\theta} - t^2 +t^{1+\theta} + t) ~|~ t \in \K \}.
\end{equation}
This has to coincide with a certain ordinary line $C_{(1,a',a'')}= \{(1+s,a'+s,a'' + (a'-1)s -s^{2+\theta}) ~|~
s \in \K \}$ (because the parallel class is preserved) with $a',a'' \in \K$. This yields the following system
of equalities~:
\begin{equation}
\left\{ \begin{array}{ll}
         t-1 =1 +s, \\
      t+1 = a'+s, \\
      -t^{2+\theta} - t^2 +t^{1+\theta} + t =   a'' + (a'-1)s -s^{2+\theta}, \end{array} \right.
\end{equation}
which simplifies to~:
\begin{equation}
\left\{ \begin{array}{ll}
         s = t+1, \\
      a' = 0 ,\\
      t =  a'' + 1 - t^\theta. \end{array} \right.
\end{equation}
If $t=0$ the last equation gives us $a'' = -1$, but if we use $t=1$, we obtain $a'' = 1$, which is a
contradiction as $a''$ is a constant. It follows that the hypothesis of one orbit is false. \qed

The following corollary follows directly~:
\begin{cor} \label{cor:gnarlcircles}
The gnarls of circles are mapped onto gnarls of circles.
\end{cor}

Using the above and Lemma~\ref{cor:sfeeruitcirkels}, one can reconstruct the spheres, giving the following
result (which is part of Main Result 2)~:
\begin{cor}
The automorphism group of $\mathcal{G_C}$ is equal to that of $\mathcal{G}$.
\end{cor}

\section{Absolute points and lines of polarities in the Ree hexagon}

We now show our Main Corollary in the formulation below. We note that our proof will not use the full strength
of our results proved so far. Indeed, we will only use Corollary~\ref{cor:gnarlcircles}. The last few lines of
the proof of the next corollary can be deleted if we use Main Result~1.

\begin{cor}
If a collineation $\sigma$  of a Moufang hexagon stabilizes the set of all absolute points of some polarity,
then it stabilizes the set of all absolute lines as well.
\end{cor}
\proof
By Theorem 7.3.4 and Theorem 7.7.2 of \cite{Mal:98}, any polarity $\rho$ of a Moufang hexagon is associated to
a Ree group, so it is a polarity of the associated Ree hexagon.

As mentioned before, a circle $C$ of the Ree geometry is the set of absolute points at distance 3 from a line
$M$, not going through an absolute point. The collineation $\sigma$ maps this set to the set of absolute points
at distance 3 from $M^\sigma$, which is again a circle since $M^\sigma$ clearly is not incident with any
absolute point (as $\sigma$ stabilizes the set of absolute points). It follows that $\sigma$ induces an
automorphism of $\mathcal{G_C}$. The gnarl of $C$ is the absolute point $x$ such that the corresponding
absolute line $x^\rho$ intersects $M$. Corollary~\ref{cor:gnarlcircles} now implies that the absolute line
$(x^\sigma)^\rho$ intersects $M^\sigma$. As $(x^\rho)^\sigma$ also contains $x^\sigma$ and intersects
$M^\sigma$, it follows that $(x^\sigma)^\rho= (x^\rho)^\sigma$. This means that the absolute line $x^\rho$ is
mapped to another absolute line. Varying $C$ we now see that the set of all absolute lines is stablized by
$\sigma$. \qed

\section{Automorphism group of the truncated Ree geometry restricted to points and spheres}\label{section:sferen}
\textbf{General idea.} | Let $\mathcal{G_S}=(\cP,\cS)$ be the truncated Ree geometry with $\cS$ the set of
spheres. We again prove that the gnarls of spheres have to be mapped onto gnarls of spheres. As a consequence
one can recognize certain automorphisms of the Ree geometry generating the Ree group, from which the circles
can be reconstructed giving us the full Ree geometry $\mathcal{G}$ and its automorphism group.

We denote by $\mathcal{G_S}'$ the derived geometry in $(\infty)$ (so the point set is
$\cP\setminus\{(\infty)\}$ and the blocks are the vertical and ordinary planes, as defined in
Section~\ref{derived}).

We start with some small observations~:
\begin{lemma}\label{lemma:intersect}
A vertical plane and an ordinary plane always intersect.
\end{lemma}
\proof By transitivity we can suppose that the vertical plane is given by $$P_a = \{(a,t',t'')~|~ t',t'' \in
\K\}, \mbox{ with }a \in \K$$ and the ordinary plane by $S_{(0,0,0)}$, which is the set $$\left\{\left(\frac{x''^\theta
-x'x''}{x''^2+x'^{1+\theta}},\frac{-x'^\theta}{x''^2+x'^{1+\theta}},\frac{-x''}{x''^2+x'^{1+\theta}}\right)~|~\K\times\K\ni (x',x'')\ne(0,0)\right\}\cup\{(0,0,0)\}.$$ If $a=0$, then $(0,0,0)\in P_a\cap S_{(0,0,0)}$. If $a\neq 0$,
then, putting $x' =0$ and $x'' = a^{-2-\theta}$ in the formula of $S_{(0,0,0)}$ gives the point
$(a,0,-a^{2+\theta})$, which is also a point of $P_a$. \qed

\begin{lemma}\label{lemma:intersect3}
The intersection of $P_0$ and $S_{(0,0,0)}$ is given by the set $\{(0,t,0)~|~t \in K \} \cup \{(0,t^{\theta
-1},t)~|~ t \in \K \backslash \{0\}\}$.
\end{lemma}
\proof Using the representations of $P_0 = \{(0,t,t')~|~t,t' \in \K \}$ and $S_{(0,0,0)} =$
$$\left\{\left(\frac{x''^\theta -x'x''}{x''^2+x'^{1+\theta}},\frac{-x'^\theta}{x''^2+x'^{1+\theta}}
,\frac{-x''}{x''^2+x'^{1+\theta}}\right)~|~\K\times \K\ni (x',x'')\ne(0,0)\right\}\cup\{(0,0,0)\},$$ we see that the points
of the intersection are determined by the equation $x''^\theta -x'x''=0$. The solutions of this equation are
given by $x''=0$ or $x' =x''^{\theta-1}$. The first set of solutions gives us $\{(0,t,0)|t \in K \}$, the
second $\{(0,t^{\theta -1},t)| t \in \K \backslash \{0\}\}$. \qed

Note that $P_0$ is the disjoint union of affine unital blocks. Indeed, the affine blocks
$\{(0,t,b)~|~t\in\K\}$, with $b\in\K$, partition $P_0$. It is now clear that the intersection of $S_{(0,0,0)}$
and $P_0$ contains exactly one affine unital block, and all other affine unital blocks in $P_0$ share exactly
one point with that intersection.

\begin{lemma}\label{lemma:intersect2}
The ordinary planes $S_{(0,0,0)}$ and $S_{(0,a',a'')}$ with $a',a'' \in \K$ intersect.
\end{lemma}
\proof Since $(0,a',a'') \in U_\infty$ maps $P_0$ to itself and $S_{(0,0,0)}$ to $S_{(0,a',a'')}$, it follows
from the paragraph preceding this lemma that $P_0\cap S_{(0,a',a'')}$ contains an affine unital block $B$. But
from that same paragraph also follows that $B$ shares a point with $P_0\cap S_{(0,0,0)}$. That point is hence
contained in $S_{(0,0,0)}\cap S_{(0,a',a'')}$. \qed

The above lemmas now allow us to prove the following analogue to Lemma~\ref{gnarlscircles}.

\begin{lemma}\label{lemma:gnarlsfeer}
The full group $G$ of automorphisms of $\mathcal{G_S}'$ has two orbits on the planes, which are the vertical
and the ordinary planes.
\end{lemma}
\proof As with the case of points and circles, it suffices to prove that the planes can not be all in one
orbit. So suppose this is the case.

We call two vertical or ordinary planes \emph{parallel} if they are disjoint or equal. By the transitivity
assumption on the planes and Lemma~\ref{lemma:intersect}, for each point $p$ (different from $(\infty)$) and
plane $P$, there is exactly one plane $Q$ parallel to $P$ and containing $p$. Let $\varpi$ be the parallel
class where $S_{(0,0,0)}$ belongs to. Because $U_\infty$ preserves parallelism and acts regularly on the
ordinary planes, the stabilizer $V$ of $\varpi$ in $U_\infty$ acts regularly on the planes in $\varpi$ and
$S_{(a,a',a'')} \in \varpi$ if and only if $(a,a',a'') \in V$.

Let $g=(a,a',a'') \in U_\infty$ be a non-trivial element of $V$. Then, in view of Lemma~\ref{lemma:intersect2},
$a$ has to be different from $0$. But as $V$ is a group, $g^3 = (0,0,-a^{2+\theta} )$ is also a non-trivial
element of $V$, which does have as first coordinate 0, so the hypothesis is false. \qed

\begin{lemma}
In $\mathcal{G_S}'$ the affine unital blocks are (geometric) invariants.
\end{lemma}
\proof We will denote the intersection of a vertical plane through the point $p$ with the ordinary plane with
gnarl $p$ by $W_{p}$. The sets $W_p$ are invariants of the geometry by virtue of Lemma~\ref{lemma:gnarlsfeer}.
Lemma~\ref{lemma:intersect3} implies that the affine unital block through $p$ is contained in $W_{p}$.

By transitivity,  it suffices to construct the affine unital block $B$ through $(0,0,0)$. Let $p \in
W_{(0,0,0)}$ be a point different from $(0,0,0)$. If $p$ lies on $B$, then $W_{(0,0,0)} \cap W_{p}$ contains
$B$ itself and thus at least $4$ points (as $|\K| > 3$). Now suppose $p\notin B$, so $p = (0,k^{\theta - 1},k)$
for a certain $k \in \K$ different from 0. Using $(0,k^{\theta - 1},k)$ as an element of $U_\infty$ and
Lemma~\ref{lemma:intersect3} we calculate that $W_{p} = \{(0,t+k^{\theta - 1},k)~|~t \in \K \} \cup
\{(0,t^{\theta -1}+k^{\theta - 1} ,t+k)~|~ t \in \K \backslash \{0\}\}$. The intersection $W_{(0,0,0)} \cap
W_{p}$ contains two obvious intersection points on the affine unital blocks contained in either $W_{(0,0,0)}$
and $W_{p}$. To look for more intersection points we need to investigate whether or not  it is possible to have
$(0,t^{\theta -1}+k^{\theta - 1} ,t+k) = (0,s^{\theta -1},s)$ for certain $s,t \in \K \backslash \{0\}$.
Equality on the third coordinate gives us $t+k=s$, the second gives us~:
\begin{align}
s^{\theta-1} = t^{\theta-1} +k^{\theta-1} &\Leftrightarrow (t+k)^{\theta-1} = t^{\theta-1} +k^{\theta-1} \\
&\Leftrightarrow t^{2-\theta} = -k^{2-\theta}.
\end{align}
If we raise both hand sides of the last equation to the power $2+\theta$, then we obtain $t=-k$, implying
$s=0$, a contradiction.

Thus in this case we have that $|W_{(0,0,0)} \cap W_{p}| = 2$. This allows us to recognize the points of the
affine unital block through $(0,0,0)$ as those for which  $|W_{(0,0,0)} \cap W_{p}| > 2$. \qed

\begin{lemma}
In $\mathcal{G_S}$, the circles of $\cG$ are invariants.
\end{lemma}
\proof Let $p$ and $q$ be two different points of $\mathcal{G_S}$ and put $G$ equal to the full automorphism
group of $\mathcal{G_S}$. Then we first want to determine the elements of $G$ which fix $p$ and all the blocks
of the unital through $p$ within the sphere with gnarl $p$ through $q$. We will denote this group by
$G_{[p,q[}$.

By 2-transitivity we can suppose that $p=(\infty)$ and $q=(0,0,0)$. The aim is to prove that
$G_{[(\infty),(0,0,0)[} = \{(0,t,0)~|~t \in \K\}=:H$. It is easy to see that these automorphisms satisfy the
needed properties and act transitively (even regularly) on the points of the affine unital block $B$ through
$(0,0,0)$. Suppose there is another automorphism $g$ which satisfies the properties. Then, possibly by
composing with a suitable element of $H$, we may assume that $g$ fixes $(0,0,0)$. This implies that the sphere
with gnarl $(0,0,0)$ through $(\infty)$ is also fixed. By lemma~\ref{lemma:intersect3} the points
$(0,k^{\theta-1},k)$ with $k \in \K \backslash \{0\}$ are also fixed, thus also the blocks through $(0,0,0)$ in
the sphere with gnarl $(0,0,0)$ through $(\infty)$, which makes the situation symmetric in both points. We can
also let the fixed points of the form $(0,k^{\theta-1},k)$ play the role of $(0,0,0)$, which yields the fixed points
$(0,k_1^{\theta-1} + k_2^{\theta-1} + \dots + k_n^{\theta-1}, k_1 +k_2 +\dots+k_n )$ with $k_i \in \K
\backslash \{0\}$ by repeating the argument. Choosing $n=3$ and $k_1 = -k_2 = k_3 =k$ with $k \in \K \backslash
\{0\}$ gives us the fixed points $(0,0,k)$ for all $k \in \K$.

Interchanging the roles of $(\infty)$ and $(0,0,0)$, we get the fixed points $(k,0,-k^{2+\theta})$ (to calculate
these observe that $(0,0,k)$ are the points different from $(\infty)$ on the circle with gnarl $(\infty)$
through $(0,0,0)$, interchanging gives us the points different from $(0,0,0)$ on the circle with gnarl
$(0,0,0)$ through $(\infty)$). If we let a fixed point $(0,0,l)$ with $l \in \K$ play the role of $(0,0,0)$, we
obtain that all the points of the form $(k,0,l)$ with $k,l \in K$ are fixed points. On each affine unital block
lies a point of this form, so all affine unital blocks are fixed, and by symmetry also the blocks of the Ree
unital through $(0,0,0)$. It follows that all points are fixed points, and that $g$ is the identity.

The above proves that $G_{[p,q[}$ is a subgroup of the root group $U_p$ and hence, if $|\K|>3$, also a subgroup
of the simple Ree group $\mathsf{R}'(\K,\theta)$. The group $K$ generated by all groups of the form $G_{[p,q[}$ is a normal subgroup of
this Ree group (indeed, if $g$ is a automorphism of $\mathcal{G_S}$, then $G_{[p,q[}^g=G_{[p^g,q^g[}$). By
simplicity, $K$ coincides with $\mathsf{R}'(\K,\theta)$. Now, by~\cite{haotea}, the root groups of $K$ are the unique
unipotent subgroups of $K$. Hence we can recover these root groups and consequently also the circles
constructed from these root groups.

If $|\K|=3$, then $K$ is a normal subgroup of the Ree group $\mathsf{R}(3)$ over the field with $3$ elements. But the
groups $G_{[p,q[}$ do not belong to the simple Ree group. Hence, it is easy to see that $K$ coincides with the
Ree group $\mathsf{R}(3)$ and, as above, we can again recover the circles. \qed

We have proved~:
\begin{cor}
The automorphism group of $\mathcal{G_S}$ coincides with that of $\mathcal{G}$.
\end{cor}

This completes the proof of Main Result 2.

Address of the Authors:

Ghent University\\
Department of Pure Mathematics and Computer Algebra\\
Krijgslaan 281, S22,\\
B-9000 Gent\\
Belgium.

\texttt{fhaot@cage.UGent.be, kstruyve@cage.UGent.be, hvm@cage.UGent.be}

\end{document}